\documentclass[final]{elsarticle}

\usepackage{amssymb,latexsym}
\usepackage{amsmath}
\usepackage{maple2e}
\usepackage{color}

\journal{: \; Expositiones Mathematicae}

\begin{document}

\begin{frontmatter}

\title{Approximate expressions for mathematical constants from PSLQ algorithm: a simple approach and a case study}  

\author{F. M. S. Lima}

\address{Institute of Physics, University of Brasilia, P.O. Box 04455, 70919-970, Brasilia-DF, Brazil}


\ead{fabio@fis.unb.br}


\begin{abstract}
In this note, I present a simple PSLQ code for finding null linear combinations, with the best rational coefficients, of mathematical constants, within some prescribed precision. As an example, I explore approximate expressions for the Ap\'{e}ry's constant $\,\zeta{(3)} = \sum_{n\ge1}{\,1/n^3}$, an irrational number to which no exact, finite closed-form expression is known. 
On taking into account a suitable search basis, I have found a simple expression for $\,\zeta{(3)}\,$ accurate to $21$ decimal places, which is triply more accurate than the best previous one. As the short \emph{Maple}$^\mathrm{TM}$ code presented here can be easily adapted to study other constants, I decided to supply it to encourage the readers to conduct their own computational experiments, as well as to adopt it in projects of numerical analysis, number theory, or linear algebra.
\newline
\end{abstract}


\end{frontmatter}

\begin{quotation}
\noindent ``{\it In fact, numerical experimentation is crucial to Number Theory, perhaps more so than to other areas of mathematics. (...)  Indeed, as Cassels has said, to a large degree Number Theory is an experimental science.''}

\begin{flushright} F. R. Villegas \quad \end{flushright}
\end{quotation}

\section{Introduction}

The Ap\'{e}ry's constant is defined as the real number to which the infinite series $\,\sum_{n=1}^\infty{\frac{1}{n^3}}$ converges (i.e., $1.2020569\ldots$) and it is so designated in honor to R. Ap\'{e}ry, who proved in 1978 that this number is irrational~\cite{Apery}. The convergence of this series, though slow,\footnote{The partial sum $\,\sum_{n=1}^{N} \, 1/n^3$, $N$ being any positive integer, of course yields a rational approximation to $\,\zeta{(3)}$. However, even for $N$ as large as $100$ yields only \emph{four} correct decimal places.} is guaranteed by the Cauchy's integral test, a result that remains valid for $\,\sum_{n=1}^\infty{{\,1/{n^s}}}$ for any complex number $s$ with $\,\Re{(s)}>1$, a broader domain in which this series is defined as $\,\zeta(s)$, the Riemann's zeta function. 
Ap\'{e}ry's constant can then be identified with $\zeta(3)$.

For positive integer values of $s$, $s>1$, Euler was the first to derive (in 1734) an exact closed-form expression for $\zeta(s)$, namely $\,\zeta(2) = \pi^2/6$, the solution of the Basel problem (see Ref.~\cite{Ayoub} and references therein).  Some years later (the result was found in 1739, but published only in 1750), he succeeded in extending his result to all \emph{even} values of $s$~\cite{Apostol}:
\begin{eqnarray*}
\zeta{(2\,n)} = (-1)^{n-1} \, \frac{2^{2 n-1} \,B_{2 n}}{(2 n)!} \; \pi^{2 n} ,
\end{eqnarray*}
where $n$ is a positive integer and $\,B_{2 n} \in \mathbb{Q}$ are Bernoulli numbers.
For \emph{odd} values of $s$, however, no exact closed-form expression is currently known.\footnote{In a recent work, it is claimed that $\,\zeta{(2\,n +1)}\,$ is not a rational multiple of $\,\pi^{2n+1}$, $\,n\,$ being a positive integer (see Theorem~22 of Ref.~\cite{chines}).}

The increase of interest in $\zeta{(3)}$, which comes from both pure and applied mathematics,\footnote{For instance, given three integers chosen at random, the probability that no common factor will divide them all is ${\,1/\zeta(3)}$.  Also, if $n$ is a power of $2$, then the number $\#(n)$ of distinct solutions for $n=p+x\,y$ with $p$ prime and $x$, $y$ positive integers obeys the asymptotic relation ${\,\#(n)/n} \sim 105 \, {\,\zeta(3)/\left(2\,\pi^4\right)}$.  It also arises in a number of physical problems, including the computation of the electron's gyromagnetic ratio.} has stimulated its high-precision numerical computation~\cite{Prevost}, 
as well as the search for \emph{simple} approximate expressions~\cite{Zagier}. Let us adopt a reasonable criterium for the adjective ``simple,'' in the context of finite approximate expressions.  Here, it will designate closed-form expressions containing \emph{a few} terms/factors composed by other known mathematical constants and \emph{a few} integer numbers with \emph{small} absolute value. This criterium is, of course, vague due to the forms ``a few'' and ``small''.  To make it \emph{not-so-imprecise}, ``a few'' will mean less than, say $\,10$, and ``small'' will mean less than, say $\,500$.  A such approximation was presented by Galliani (2002), namely~\cite{Wolfram}
\begin{equation}
\zeta(3) \approx \frac{1}{\sqrt[3]{\gamma}} + \frac{1}{\pi^4} \, \left(1+2\,\gamma-\frac{2}{130+\pi^2}\right)^{\!-3} ,
\label{eq:Galli}
\end{equation}
where $\,\gamma\,$ is the Euler-Mascheroni constant, which is accurate to $4$ decimal places.  Another nice simple approximation is
\begin{equation}
\zeta(3) \approx \sqrt[4]{\, \gamma + \frac{71}{47}\,} \: ,
\label{eq:bela}
\end{equation}
due to Hudson (2004), which is accurate to $7$ places~\cite{Wolfram}.  Among the many approximations for $\zeta{(3)}$ presented by Hudson, the most accurate is~\cite{Wolfram}
\begin{equation}
\zeta(3) \approx 525587^{\; 1/\sqrt{5123}} \, ,
\label{eq:boa}
\end{equation}
which yields $12$ correct decimal places and, clearly, is not a \emph{simple} approximate expression (in our terminology).  The same for $\zeta(3) \approx \frac{97525}{2515594} \, \pi^3$, which I have found on searching for a direct integer relation between $\zeta(3)$ and $\pi^3$. 
On trying to reach greater accuracy, however, one soon observes that \emph{it is very difficult to avoid the appearance of large integers}.
\medskip

In this note, I show how to use PSLQ algorithm for finding the best rational coefficients, within a desired precision, of a null linear combination of a finite number of real constants. As an example, I explore approximate expressions for the Ap\'{e}ry's constant $\,\zeta{(3)}$, a number to which no closed-form expression in terms of a finite combination of elementary functions of known constants is known. For this, I choose a suitable search basis composed by numbers which seem to be closely related to $\,\zeta{(3)}$, namely $\,\pi$, $\,\ln{2}\,$, $\,\ln{(1+\sqrt{2}\,)}$, and $G$ (the Catalan's constant). This yields a simple expression for $\,\zeta{(3)}\,$ accurate to $21$ decimal places. The short \emph{Maple}$^\mathrm{TM}$ code used in this computation is included to stimulate the readers to conduct their own experiments.

\section{Searching for integer relations: the PSLQ algorithm}

An important task in experimental mathematics is to search for integer relations involving a finite set of known numbers.  An integer relation algorithm is a computational scheme that, for a given real vector $\mathbf{x} = \left(x_1,x_2,\ldots,x_n\right)$, $n>1$, it either finds a nonnull vector of integers $\mathbf{a} = \left(a_1,a_2,\ldots,a_n\right)$ such that  $a_1 \, x_1 + a_2 \, x_2 + \ldots + a_n \, x_n = 0$ or else establishes that there is no such integer vector within a ball of some radius about the origin.\footnote{The metric is the Euclidean norm $\,\sqrt{a_1^2 +a_2^2 + \ldots +a_n^2}\:$.}

Presently, the best algorithm for detecting integer relations is the PSLQ  algorithm (acronym for \emph{Partial Sum of Least sQuares}) introduced by Ferguson and Bailey (1992)~\cite{Ferguson}. A simplified formulation for this algorithm, equivalent to the original one, was subsequently developed by Ferguson and co-workers (1999)~\cite{Bailey99}. This more efficient version is currently implemented in both \emph{Maple}$^\mathrm{TM}$ and \emph{Mathematica}$^\mathrm{TM}$, two of the most popular mathematical softwares. This version of PSLQ, optimized with certain reduction schemes,
was named one of the `ten algorithms of the century' in Ref.~\cite{ALGcentury}. 

In short, PSLQ operates as follows. Given a vector $\mathbf{x}$ of $n$ given real numbers, input as a list of floating-point (FP) numbers, the algorithm uses \emph{QR decomposition} in order to construct a series of matrices $A_m$ such that the absolute values of the entries of the vector $\mathbf{y}_m = A_m^{-1} \cdot \mathbf{x}$ decrease monotonically. At any given iteration, the largest and smallest entries of $\mathbf{y}_m$ usually differ by no more than a few orders of magnitude. When the desired integer relation is detected, the smallest entry of $\mathbf{y}_m$ abruptly decreases to roughly the computer working precision $\epsilon$ and the relation is given by the corresponding column of $A_m^{-1}$.  This numerically stable matrix reduction procedure, together with some techniques that allow machine arithmetic to be used in many intermediate steps, usually yields a rapid convergence, which makes PSLQ faster than other concurrent algorithms~\cite{bertok}. In fact, if the elements of $\mathbf{x}$ are linearly dependent over $\mathbb{Q}$, then PSLQ will find an exact integer relation between them, for a sufficiently precise input. For most applications, high-precision arithmetic is required, which stems from the fact that if one wishes to recover a relation involving $n$ known real numbers, with coefficients accurate to $d$ digits, then the input vector $\mathbf{x}$ must be specified to at least $n \times d$ digits and one must employ FP arithmetic accurate to at least $n \times d$ digits, too.
When a relation is detected, the ratio between the smallest and the largest entry of the vector $A^{-1} \cdot \mathbf{x}$ can be taken as a ``confidence level'' that the relation is true (i.e., exact) and not an artifact of insufficient numerical precision. Very small ratios at detection certainly indicate the result is probably true.\footnote{In addition to possessing good numerical stability, PSLQ is guaranteed to find an integer relation in a number of iterations bounded by a polynomial in $n$.}


\section{Closed-form expressions via PSLQ}

Since an efficient PSLQ routine is available as part of a \emph{Maple}$^\mathrm{TM}$ package named \verb"IntegerRelations", then simple short codes can be written in this high-level language for finding approximate expressions for mathematical constants. For illustrating this, let me list the source code I have written for finding an approximate expression for $\zeta(3)$.

{\small
\def\emptyline{\vspace{12pt}}
\DefineParaStyle{Maple Output}
\DefineCharStyle{2D Math}
\DefineCharStyle{2D Output}



\begin{maplegroup}
\begin{mapleinput}
\mapleinline{active}{1d}{Digits := 24:  # The number of digits for FP numbers}{%
}
\end{mapleinput}
\end{maplegroup}

\begin{maplegroup}
\begin{mapleinput}
\mapleinline{active}{1d}{with(IntegerRelations): # Call the package containing PSLQ}{%
}
\end{mapleinput}
\end{maplegroup}

\begin{maplegroup}
\begin{mapleinput}
\mapleinline{active}{1d}{xSymb:=[Zeta(3),1,Pi^2*ln(2),Pi*ln(2)^2,ln(2)^3,ln(1+sqrt(2))^3,Pi*Catalan];}{%
}
\end{mapleinput}
\mapleresult
\begin{maplelatex}
\mapleinline{inert}{2d}{xSymb := [Zeta(3), 1, Pi^2*ln(2), Pi*ln(2)^2, ln(2)^3, ln(1+2^(1/2))^3, Pi*Catalan];}{%
\[ \mathit{xSymb} := [\zeta (3), 1, \,\pi ^{2}\,\mathrm{ln}(2), \,\pi \,\mathrm{ln}(2)^{2}, \,\mathrm{ln}(2)^{3}, \,\mathrm{ln}(1 + \sqrt{2})^3, \,\pi \,\mathit{Catalan}] \] }
\end{maplelatex}
\end{maplegroup}

\begin{maplegroup}
\begin{mapleinput}
\mapleinline{active}{1d}{n := nops(xSymb);  # The number of elements in xSymb}{%
}
\end{mapleinput}
\mapleresult
\begin{maplelatex}
\mapleinline{inert}{2d}{n := 7;}{%
\[ n := 7 \] }
\end{maplelatex}
\end{maplegroup}

\begin{maplegroup}
\begin{mapleinput}
\mapleinline{active}{1d}{x:=evalf(xSymb):  # Convert to FP numbers}{%
}
\end{mapleinput}
\end{maplegroup}

\begin{maplegroup}
\begin{mapleinput}
\mapleinline{active}{1d}{a:=PSLQ(x);  # Applies PSLQ algorithm to vector x}{%
}
\end{mapleinput}
\mapleresult
\begin{maplelatex}
\mapleinline{inert}{2d}{a := [10, 394, -11, 283, -472, -209, -186];}{%
\[ a := [10, \,394, \,-11, \,283, \,-472, \,-209, \,-186] \] }
\end{maplelatex}
\end{maplegroup}

\begin{maplegroup}
\begin{mapleinput}
\mapleinline{active}{1d}{soma:=0:}{%
}
\end{mapleinput}
\end{maplegroup}

\begin{maplegroup}
\begin{mapleinput}
\mapleinline{active}{1d}{for k from 1 to n do \;  soma:=soma+a[k]*xSymb[k]; \:  end do: }{%
}
\end{mapleinput}

\end{maplegroup}

\begin{maplegroup}
\begin{mapleinput}
\mapleinline{active}{1d}{aprox:=solve(soma=0, Zeta(3));}{%
}
\end{mapleinput}
\mapleresult
\begin{maplelatex}
\mapleinline{inert}{2d}{aprox:=-5/197+11/394*Pi^2*ln(2)-283/394*Pi*ln(2)^2+236/197*ln(2)^3+209/394*ln(1+2^(1/2))^3 +93/197*Pi*Catalan;}{%
\[ \mathit{aprox} :=  - {\displaystyle \frac {5}{197}}  +
{\displaystyle \frac {11}{394}} \,\pi ^{2}\,\mathrm{ln}(2) -
{\displaystyle \frac {283}{394}} \,\pi \,\mathrm{ln}(2)^{2} +
{\displaystyle \frac {236}{197}} \,\mathrm{ln}(2)^{3} +
{\displaystyle \frac {209}{394}} \,\mathrm{ln}(1 + \sqrt{2})^3 \\
+ {\displaystyle \frac{93\,\pi \,\mathit{Catalan}}{197}}
\] }
\end{maplelatex}
\end{maplegroup}

\begin{maplegroup}
\begin{mapleinput}
\mapleinline{active}{1d}{evalf( aprox, 30 );  # Our approximation}{%
}
\end{mapleinput}
\mapleresult
\begin{maplelatex}
\mapleinline{inert}{2d}{1.20205690315959428539958993430;}{%
\[ {\color{blue} 1.202056903159594285399}58993430 \] }
\end{maplelatex}
\end{maplegroup}

\begin{maplegroup}
\begin{mapleinput}
\mapleinline{active}{1d}{evalf( Zeta(3), 30 );  # Exact value (for comparison)}{%
}
\end{mapleinput}
\mapleresult
\begin{maplelatex}
\mapleinline{inert}{2d}{1.20205690315959428539973816151;}{%
\[ {\color{blue} 1.20205690315959428539973816151} \] }
\end{maplelatex}
\end{maplegroup}
}


\bigskip

With this \emph{Maple}$^\mathrm{TM}$ routine,  I have found the following simple approximate expression for the Ap\'{e}ry's constant:
\begin{equation}
\zeta(3) \approx -\frac{5}{197} \,+\frac{11}{394} \, \pi^2 \ln{2} \,-\frac{283}{394} \, \pi \, \ln^2{2} \,+\frac{236}{197} \, \ln^3{2} \,+\frac{209}{394} \, \alpha^3 +\frac{93}{197} \, \pi \, G , \;
\label{eq:my21}
\end{equation}
where $\alpha \equiv \ln{(1 +\sqrt{2}\,)}$, which is accurate to $21$ digits, as the reader can check by comparing the last two numerical outputs above (the blue digits are correct).  This approximate expression, Eq.~(\ref{eq:my21}), is much more accurate than the best previous expressions and this is why it has been included in the popular webpage \emph{Wolfram MathWorld}~\cite{Wolfram}.
\newline

All that rests is to explain the motivation for choosing the vector
\begin{eqnarray*}
\mathbf{x} = \left( 1, \, \pi^2 \, \ln{2}, \, \pi \, \ln^2{2}, \, \ln^3{2}, \, \ln^3{( 1+\sqrt{2} \,)}, \, \pi\,G \right)
\end{eqnarray*}
as the search basis for $\zeta(3)$.  This comes primarily from a conjecture by Euler (1785) that $\,\zeta(3) = a \, \ln^3{2} + b \, \pi^2 \,\ln{2}\,$, for some $\, a,b \in \mathbb{Q}$~\cite{EulerLog2}. After some experiments with a code similar to the above one, I could not find any such pair of rational coefficients for composing an \emph{exact} closed-form expression for $\zeta(3)$. Not even a \emph{simple} expression was found with these two terms only.\footnote{I am currently testing a conjecture by Connon (2008) that either $a$ or $b$ could contain a factor of $\sqrt{2}$, or another small surd, or perhaps $\,\ln{(2 \pi)}$~\cite{Connon}.}  I was then inclined to improve the basis by including other \emph{weight-3} constants in virtue of their appearance in some exact results involving $\zeta(3)$, e.g.:\footnote{The word \emph{weight}, here, follows the definition introduced by Boros (see Ref.~\cite{BookInts}, p.~203).}

\begin{description}

  \item (i) {Special values of non-elementary functions:

   $\mathrm{Li}_3\left(\frac12 \right) = \frac78 \,\zeta(3) -\frac{1}{12}\,\pi^2 \, \ln{2}  +\frac16\, \ln^3{2}\,$, $\;\Re{\left(\mathrm{Li}_3\left(\frac{1+i}{2}\right)\right)} = \frac{35}{64} \,\zeta(3) +\frac{1}{48} \, \ln^3{2} -\frac{5}{192} \, \pi^2 \, \ln{2}$~\cite{brod}, $\;\pi^2\, \psi^{(-4)}(1) = \frac12 \pi^2 \ln{A} +\frac{1}{12}\,\pi^2\,\ln{2} +\frac{1}{12}\,\pi^2\,\ln{\pi} +\frac18 \, \zeta(3)$~\cite{WolfZ3} (typo corrected);}

  \item (ii) {Infinite series:

   $\sum_{k=1}^{\infty}{\frac{(-1)^{k+1}}{k^3\,2^k\,\binom{2k}{k}}} = \frac14 \zeta(3) -\frac16 \, \ln^3{2}$~\cite{Huvent}, $\;\pi^2 \sum_{k=1}^{\infty}{\frac{\zeta(2k)}{(k+1)\,2^{4k}}} = \frac12 \pi^2 -\frac12 \, \pi^2 \,\ln{2} +\frac{35}{4}\,\zeta(3) -4 \pi\,G$~\cite{Cho}, $\;\pi^4 \sum_{m=0}^{\infty}{\frac{(-1)^{m}\, \pi^{2m}\,E_{2m+1}(1)}{(2m+5)!}} = \frac13 \pi^2 \ln{2} {-\frac32 \zeta(3)}$~\cite{Dancs};}

  \item (iii) {Definite integrals:

   $\;\int_{0}^{\pi/4}{x\,\ln{(\cos{x})} \, dx} = -\frac{21}{128} \, \zeta(3) +\frac18 \, \pi\,G -\frac{1}{32}\,\pi^2\,\ln{2}$~\cite{meuInt}, $\;\int_{0}^{\pi/4}{x^2\,\tan{x} \, dx} = -\frac{21}{64} \, \zeta(3) +\frac14 \, \pi\,G -\frac{1}{32}\,\pi^2\,\ln{2}$~\cite{meuInt}, $\;\int_{0}^{\pi/2}{\,x^2/\sin{x} \, dx} = 2\,\pi\,G -\frac72 \, \zeta(3)$~\cite{Cho}, $\;\int_{0}^{1}{ \frac{(\arcsin{t})^2}{t} \, dt} = \frac14 \pi^2 \ln{2} -\frac78 \,\zeta(3)$ (Eq.(6.6.25) in Ref.~\cite{BookInts}).} 
\end{description}
Here, $\mathrm{Li}_3(z) := \sum_{n=1}^{\infty} {\,z^n/n^3}$ is the trilogarithm function, $\psi^{(-4)}(z)$ is the polygamma function (extended to negative indexes, according to Ref.~\cite{Espinosa}), and $\,E_{2m+1}(x)\,$ is the Euler polynomial of degree $\,2m+1$.

With respect to the term with $\,\ln^3{(1+\sqrt{2}\,)}$, which is the logarithm of a non-null algebraic number distinct from unit, thus a \emph{weight-3} constant, I have included it in my search basis because the number $\,\mathrm{arcsinh}(1) = \ln{(1+\sqrt{2} \,)}\,$ emerges in the coordinates of the vertices of a cusped hyperbolic cube whose volume is $\,\frac78 \: \zeta(3)$, as I have found on exploring certain triple integrals over the unit cube $[0,1]^3$, as suggested in the end of Ref.~\cite{meuHyp}.


\end{document}